\documentclass[a4paper,12pt]{article}

\usepackage{amsfonts}
\usepackage{amscd,color}
\usepackage{amsmath,amsfonts,amssymb,amscd}
\usepackage{mathrsfs}
\usepackage{indentfirst,graphicx,epsfig}
\usepackage{graphicx}
\input{epsf}
\usepackage{epstopdf}
\usepackage{caption}
\usepackage{ifpdf}
\usepackage{fullpage}
\usepackage{epsfig}
\usepackage{latexsym}
\usepackage{caption}
\usepackage{amsmath,amsthm}
\usepackage{amsfonts}
\usepackage{amssymb}
\usepackage{graphicx}
\usepackage{graphics}
\usepackage{bm}
\usepackage{color}
\usepackage{subfigure}
\graphicspath{{\figs}}

\makeatletter

\newcommand{\Rmnum}[1]{\expandafter\@slowromancap\romannumeral #1@}
\makeatother
\makeatletter
\let\@fnsymbol\@arabic
\makeatother

\begin{document}
\newtheorem{theorem}{Theorem}[section]
\newtheorem{observation}[theorem]{Observation}
\newtheorem{corollary}[theorem]{Corollary}
\newtheorem{algorithm}[theorem]{Algorithm}
\newtheorem{problem}[theorem]{Problem}
\newtheorem{question}[theorem]{Question}
\newtheorem{lemma}[theorem]{Lemma}
\newtheorem{proposition}[theorem]{Proposition}

\newtheorem{definition}[theorem]{Definition}
\newtheorem{guess}[theorem]{Conjecture}
\newtheorem{claim}[theorem]{Claim}
\newtheorem{example}[theorem]{Example}
\makeatletter
  \newcommand\figcaption{\def\@captype{figure}\caption}
  \newcommand\tabcaption{\def\@captype{table}\caption}
\makeatother

\newtheorem{acknowledgement}[theorem]{Acknowledgement}

\newtheorem{axiom}[theorem]{Axiom}
\newtheorem{case}[theorem]{Case}
\newtheorem{conclusion}[theorem]{Conclusion}

\newtheorem{conjecture}[theorem]{Conjecture}
\newtheorem{criterion}[theorem]{Criterion}
\newtheorem{exercise}[theorem]{Exercise}
\newtheorem{notation}[theorem]{Notation}
\newtheorem{solution}[theorem]{Solution}
\newtheorem{summary}[theorem]{Summary}
\newtheorem{fact}[theorem]{Fact}
\newtheorem{remark}[theorem]{Remark}

\newcommand{\pp}{{\it p.}}
\newcommand{\de}{\em}
\newcommand{\mad}{\rm mad}

\newcommand*{\QEDA}{\hfill\ensuremath{\blacksquare}}
\newcommand*{\QEDB}{\hfill\ensuremath{\square}}

\newcommand{\qf}{Q({\cal F},s)}
\newcommand{\qff}{Q({\cal F}',s)}
\newcommand{\qfff}{Q({\cal F}'',s)}
\newcommand{\f}{{\cal F}}
\newcommand{\ff}{{\cal F}'}
\newcommand{\fff}{{\cal F}''}
\newcommand{\fs}{{\cal F},s}
\newcommand{\g}{\gamma}
\newcommand{\wrt}{with respect to }

\title{The asymptotic value of energy for matrices with degree-distance-based entries of random graphs\footnote{Supported by NSFC No.11871034 and 11531011.}
}
\author{\small Xueliang Li$^{1}$, ~ Yiyang Li$^{2}$, ~ Zhiqian Wang$^{1}$\\
\small $^{1}$Center for Combinatorics and LPMC\\
\small Nankai University, Tianjin 300071, China\\
\small Email: lxl@nankai.edu.cn; 1522686578@qq.com\\
\small $^{2}$International Institute, China Construction Bank\\
\small Xicheng District, Beijing 100033, China\\
\small Email: liyiyang@163.com\\
}

\date{ }
\maketitle

\begin{abstract}
For a graph $G=(V, E)$ and $i, j\in V$, denote the distance between $i$ and $j$ in $G$ by $D(i, j)$ and the degrees of $i$, $j$ by $d_i$, $d_j$, respectively. Let $f(D(i, j), d_{i}, d_{j})$ be a function symmetric in $i$ and $j$. Define a matrix $W_f(G)$, called the weighted distance matrix, of $G$, with the $ij$-entry $W_f(G)(i, j)=f(D(i, j), d_{i}, d_{j})$ if $i\neq j$ and $W_f(G)(i, j)=0$ if $i=j$. In this paper, we prove that if the symmetric function $f$ satisfies that $f(D(i, j), (1+o(1))np, (1+o(1))np)=(1+o(1))f(D(i, j), np, np)$, then for almost all graphs $G_p$ in the $Erd\ddot{o}s$-$R\acute{e}nyi$ random graph model $\mathcal{G}_{n, p}$, the energy of $W_f(G_p)$ is $\{(\frac{8}{3\pi}\sqrt{p(1-p)}+o(1))\cdot|f(1, np, np)-f(2, np, np)|+o(|f(2, np, np)|)\}\cdot n^{3/2}$. As a consequence, we give the asymptotic values of energies of a variety of weighted distance matrices with function $f$ from distance-based only and mixed with degree-distance-based topological indices of chemical use. This
generalizes our former result with only degree-based weights.\\[3mm]
\end{abstract}

\section[Introduction]{Introduction}

Throughout the paper, we denote a simple graph by $G=(V, E)$ with order $|V(G)|=n$. The adjacency matrix of $G$ is denoted by $A(G)$. If $e\in E$ is an edge with two ends $i$ and $j$, we say that $e=ij\in E$. We use $d_{i}$ and $D(i, j)$ to represent the degree of a vertex $i$ and the distance between two vertices $i$ and $j$ in $G$, respectively.\\

In graph theory and its applications, especially in chemistry, matrices are popularly introduced and studied, which provides algebraic perspectives on (molecular) graphs. The most classic graph matrix is the adjacency matrix, whose entries indicate the adjacency (or number of edges) between two vertices. In practical requirements of molecular chemistry, many objects are concerned with degrees of vertices or distances between pairs of vertices in a graph, and thus many interesting matrices with entries from degrees and/or distances have been introduced. One kind of such matrices comes from degree-based topological indices of chemical use, such as the Zagreb matrix \cite{gt}, ABC-matrix \cite{est} and Harmonic matrix \cite{hkb}, which are essentially the adjacency matrix weighted by a symmetric function $f(d_{i}, d_{j})$ defined on the degrees of vertices $i$ and $j$. Another kind of such matrices comes from distance-based topological indices, such as the distance matrix \cite{dll2}, Harary matrix \cite{plav} and reverse Wiener matrix \cite{van}, which are essentially  the distance matrix with a symmetric function $f(D(i, j))$ as the $ij$-entry for $i\neq j$ in $G$. In 1994, Dobrynin and Kochetova \cite{dob} put forward a new topological index determined by the values of both distances and degrees of vertices, and recently this new type of indices become more and more popular. We refer the reader to \cite{amm, feng, mm, xu} for results on the degree-distance-based indices. From this kind of indices, it is natural to define a new kind of matrices with mixed degree-distance-based  entries, since a 2-dimensional matrix contains much more data than a single index, and its algebraic property will show more structural information of a molecular. As one can see below, it is also essentially a distance matrix, and so we call it the \emph{weighted distance matrix} of a graph $G$. The definition is given as follows.\\

\begin{definition}
Let $G=(V,E)$ be a graph. Denote by $d_i$ the degree of a vertex $i$ in $G$, and by $D(i,j)$ the distance between two vertices $i,j$ in $G$. Let $f(D(i,j),d_i,d_j)$ be a function symmetric in $d_i$ and $d_j$. The weighted distance matrix $W_f(G)$ of $G$ is defined as follows: the $ij$-entry of $W_f(G)$
\begin{displaymath}
W_f(G)(i, j)=
\begin{cases}
f(D(i, j), d_i, d_j), &\quad \mbox{$i\neq j$}\\
0,&\quad \mbox{$i=j$}
\end{cases}
\end{displaymath}
\end{definition}

\begin{remark} \label{rem1} In fact, the above adjacency matrix with only degree-based weights and the distance matrix with only distance-based weights are special cases of this new type of matrix $W_f(G)$. If we set $f(D(i, j), d_i, d_j)=0$ for $D(i, j)\geq 2$, then $W_f(G)$ is the adjacency matrix with degree-based weights by function $f(1, d_i, d_j)$. If the value of the function $f(D(i, j), d_i, d_j)$ depends only on $D(i, j)$, $W_f(G)$ is the distance matrix with $ij$-entries weighted by only distance-based function $f$.
\end{remark}

Since $f$ is a symmetric function, $W_f(G)$ is a symmetric matrix and therefore it has only real eigenvalues, denoted by $\lambda_1 \geq \lambda_2\geq\cdots \geq \lambda_n$. As usual, the energy of the matrix $W_f(G)$ is defined as follows:
\begin{displaymath}
\mathscr{E}(W_f(G))=\sum_{i=1}^{n}|\lambda_{i}|,
\end{displaymath}

The concept of energy for graphs, derived from chemistry, was first introduced by Gutman \cite{gut1} in 1978. The total $\pi$-electron energy in a conjugated hydrocarbon is given by the sum of absolute values of the eigenvalues corresponding to the molecular graph. Recently, various energy of degree and/or distance based matrices have been studied, such as the Randi\'c energy and ABC energy (see \cite{lls}), distance energy \cite{dll2} and Harary energy \cite{plav}, etc. These types of energies were introduced simply by replacing the adjacency matrix $A(G)$ with the corresponding matrix from chemical indices. In this paper we will study the energy of a more general matrix, the weighted distance matrix  $W_f(G)$.\\

It is not difficult to calculate the exact value of energy of $W_f(G)$ for a concrete graph $G$ and some simple functions $f$, just by computing the eigenvalues of $W_f(G)$. But, this is impractical when $n$ is getting large and $f$ becomes more complex. So, it will be interesting to get  the asymptotic tendency of $\mathscr{E}(W_f(G))$ as $n\rightarrow\infty$. Random graphs are very suitable objects to serve this purpose. In this paper, we will study the asymptotic value of energy of $W_f(G_p)$ for random graphs $G_p$ in the classic $Erd\ddot{o}s$-$R\acute{e}nyi$ random graph model \cite{erd}, from which we can see that the asymptotic tendency of $\mathscr{E}(W_f(G))$ as a graph $G$ is getting more and more large and dense. Recall that $\mathcal{G}_{n, p}$ consists of all graphs on $n$ vertices in which the edges are chosen independently with probability $p$, where $p\in (0, 1)$ is a constant.\\

At first, we recall some results about random matrices, details of which can be found in \cite{wig1, wig2}. In 1950s, Wigner studied the limiting spectral distribution of a type of random matrices, named as \emph{Wigner matrix}, denoted by $X=\{x_{ij}\}_{i, j=1}^{n}$, which satisfies the following conditions:
\\
(i) \ $x_{ij} \ (i\neq j)$ are \emph{i.i.d.} random variables with variance $\sigma^{2}$, and $x_{ij}=x_{ji}$;\\
(ii) \ $x_{ii}$ are \emph{i.i.d.} random variables without any moment requirement.\\

The \emph{empirical spectral distribution (ESD)} of $X$ is defined by $\Phi_{X}(x)=\frac{1}{n}\cdot\sharp\{\lambda_{i}|\lambda_{i}\leq x, i=1, 2, \ldots, n\}$. Then, the energy of matrix $X$ is $\mathscr{E}(X)=n\cdot\int|x|\,d\Phi_{X}(x)$. Wigner calculated the limiting spectral distribution (LSD for short) of $X$ and obtained his famous semi-circle law; see \cite{wig1, wig2}. Throughout this paper, the expression ``a.s." means ``with probability tending to 1 as $n$ tends to infinity".
\begin{theorem} (semi-circle law)  \label{thm1} \cite{wig1, wig2}
\begin{displaymath}
\lim_{n\rightarrow\infty}\Phi_{n^{-1/2}X}(x)=\Phi(x)\quad a.s.
\end{displaymath}
i.e., with probability 1, $\Phi_{n^{-1/2}X}(x)$ converges weakly to a distribution $\Phi(x)$ as n tends to infinity. $\Phi(x)$ has the density
\begin{displaymath}
\phi(x)=
\begin{cases}
\frac{1}{2\pi\sigma^{2}}\sqrt{4\sigma^{2}-x^{2}},&\quad |x|\leq 2\sigma,\\
0, & \mbox{otherwise}.
\end{cases}
\end{displaymath}
\end{theorem}

Du, Li and Li \cite{dll1} used the semi-circle law and studied the energy of the adjacency matrix $A_p$ for random graphs $G_p$. It is easy to see that $\bar{A}_p=A_p-p(J-I)$ is a Wigner matrix with $\sigma=\sqrt{p(1-p)}$, where $J$ is the matrix of all ones and $I$ is the unit matrix. So, the limiting distribution of eigenvalues is determined by the semi-circle law. As a corollary, they calculated the energy of $G_p$ and obtained the following asymptotic result.\\
\begin{theorem} \cite{dll1}
$$\mathscr{E}(G_{p})=(\frac{8}{3\pi}\sqrt{p(1-p)}+o(1))n^{3/2}\quad a.s.$$
\end{theorem}

Unfortunately, when we consider $W_f(G_p)$ of random graphs $G_p$, a big problem arises. This matrix is no longer a Wigner matrix since the random variables $f(d_{i}, d_{j})$, $f(D(i, j))$ or $f(D(i, j), d_{i}, d_{j})$ (${i,j}\in \{1,2,\ldots,n\}$) are not independent. This limits the use of this method, and most previous results become unavailable.
However, fortunately we can still use moment method to estimate the asymptotic value of energy for large weighted random matrices, which was once applied in Wigner's paper \cite{wig1}.

\begin{lemma} \label{lem1.5} (moment method, see \cite{van}) \\
Suppose $\{Y_{n}\}_{n=1}^{\infty}$ is a sequence of random variables and $Y$ is a fixed random variable whose distribution is uniquely determined by its moments. Suppose that all of the moments $\mathbf{E}(Y_{n}^{k}), \ \mathbf{E}(Y^{k}) \ (k=1,2,\ldots)$ exist. If\\
\begin{displaymath}
\lim\limits_{n\rightarrow\infty}\mathbf{E}(Y_{n}^{k})=\mathbf{E}(Y^{k})
\end{displaymath}
for all values of $k$, then $Y_{n}$ converges to $Y$ in distribution.
\end{lemma}

The asymptotic value of energy for random graphs with degree-based weights $f(d_{i}, d_{j})$ was studied recently by Li, Li and Song in \cite{lls} and they obtained the following result.\\

\begin{theorem}
Let $f(x, y)$ be a symmetric real function. Denote by $A_{p}(f)$ the adjacency matrix of a random graph $G_p$ weighted by a degree-based function $f(d_i,d_j)$. If the function $f$ satisfies that conditions  that $|f(d_{i}, d_{j})|\leq Cn^{m}$ for some constants $C, m > 0$, and $f((1+o(1))np, (1+o(1))np)=(1+o(1))f(np, np)$ where $p\in(0, 1)$ is any fixed and independent of $n$, then for almost all random graphs $G_p$ in $\mathcal{G}_{n, p}$,\\
\begin{displaymath}
\mathscr{E}(A_{p}(f))=|f(np, np)|(\frac{8}{3\pi}\sqrt{p(1-p)}+o(1))\cdot n^{3/2} \quad a.s.
\end{displaymath}
\end{theorem}

As we pointed out, the above matrix $A_p(f)$ is essentially the adjacency matrix $A(G_p)$. However, our new matrix $W_f(G)$ is essentially the distance matrix of $G$, which has not been studied in existing literature, yet. In this paper, we will calculate the energy of weighted distance matrix $W_f(G_p)$ of $G_p\in\mathcal{G}_{n, p}$, and obtain a result similar to that in \cite{lls}, but more general by Remark \ref{rem1}.\\

Let $f(D(i, j), d_{i}, d_{j})$ be a function symmetric in $i$ and $j$  with the property:\\
\begin{displaymath}
f(D(i, j), (1+o(1))np, (1+o(1))np)=(1+o(1))f(D(i, j), np, np). \tag{*}
\end{displaymath}
We will give the limiting value of $\mathscr{E}(W_f(G_p))$ based on Wigner's moment method. As a result, we obtain the following result.\\
\begin{theorem} \label{main}
Let $f(D(i, j), d_{i}, d_{j})$ be a symmetric function satisfying the above property $(*)$. Then for almost all graphs $G_p\in \mathcal{G}_{n, p}$,
$$
\mathscr{E}(W_f(G_p))=\{(\frac{8}{3\pi}+o(1))\cdot|f(1, np, np)-f(2, np, np)|+o(|f(2, np, np)|)\}\cdot n^{3/2}\quad a.s.
$$
That is, if $f(1, np, np)/f(2, np, np)\nrightarrow 1$,\\
$$
\mathscr{E}(W_f(G_p))=|f(1, np, np)-f(2, np, np)|(\frac{8}{3\pi}\sqrt{p(1-p)}+o(1))\cdot n^{3/2}\quad a.s.
$$
and if $f(1, np, np)/f(2, np, np)\rightarrow 1$,\\
$$
\mathscr{E}(W_f(G_p))=o(1)|f(2, np, np)|\cdot n^{3/2}\quad a.s.
$$
\end{theorem}

From the above, one can see that when the values of $f(1, np, np)$ and $f(2, np, np)$ are sufficiently approximate, the energy becomes very small.

\section[Proof of Theorem \ref{main}]{Proof of Theorem \ref{main}}

As is shown in Section 1, the matrix $W_f(G_p)$ may be rather complicated since the diameter of a graph can be very large, and also many different values of distance $D(i, j)$ for pairs of vertices of $G$ are involved. However, things are not that disappointed from the probability point of view. In fact, almost all graphs have diameter two, see \cite{boll}. So, to study the asymptotic property, it suffices to deal with graphs of diameter $2$, whose weighted distance matrix $W_f(G_p)$ consists of entries with only values $0$, $f(1, d_{i}, d_{j})$ and $f(2, d_{i}, d_{j})$.

Suppose $G_p$ is a graph with diameter 2. We write $W_f(G_p)=A_{1}+A_{2}$, where
\begin{displaymath}
A_{1}(i, j)=
\begin{cases}
f(1, d_{i}, d_{j})-f(2, d_{i}, d_{j}),&\quad \mbox{$i$ and $j$ are adjacent,}\\
0,&\quad \mbox{$i$ and $j$ are nonadjacent, or $i=j$,}
\end{cases}
\end{displaymath}
and
\begin{displaymath}
A_{2}(i, j)=
\begin{cases}
f(2, d_{i}, d_{j}),&\quad \mbox{ $i\neq j$,}\\
0,&\quad \mbox{ $i=j$.}
\end{cases}
\end{displaymath}
The two matrices will be treated separately in the following subsections.\\

\subsection{Calculation of $\mathscr{E}(A_{1})$}

We will use moment method in this subsection. For convenience, we replace $f(1, d_{i}, d_{j})-f(2, d_{i}, d_{j})$ with $F(d_{i}, d_{j})$. Moreover, we assume that $f(1, np, np)\neq f(2, np, np)$. The case that $f(1, np, np)=f(2, np, np)$ will be discussed in the next subsection.

Let $M_{k}(A)=\int x^{k}\, d\Phi_{A}(x)$. The proof of semi-circle law is based on moment method. Consider Wigner matrix $X$ with properties:
(i) \ the diagonal entries are 0; (ii) \ the off-diagonal entries are \emph{i.i.d} with mean 0 and variance $\sigma^{2}$. It was proved in \cite{bai} that\\
\begin{lemma} \label{lem2.1}
\begin{displaymath}
\mathbf{E}M_{k}(n^{-1/2}X)=
\begin{cases}
\frac{(2s)!\sigma^{k}}{s!(s+1)!}+O(n^{-1}),&\quad \mbox{k=2s}\\
O(n^{-1/2}),&\quad \mbox{k=2s+1}
\end{cases}
\end{displaymath}
and\\
\begin{displaymath}
\mathbf{Var}(M_{k}(n^{-1/2}X))=O(n^{-2}).
\end{displaymath}
\end{lemma}
That is, for any fixed integer $k>0$,\\
\begin{displaymath}
\int x^{k}\,d\Phi_{n^{-1/2}X}(x)\rightarrow\int x^{k}\,d\Phi(x)\quad a.s.
\end{displaymath}
For more details, we refer the reader to \cite{bai}.

Let $\bar{A}_{1}=\frac{A_{1}}{f(1, np, np)-f(2, np, np)}-p(J-I)$, and $\bar{F}(d_{i_{1}}, d_{i_{2}})=\frac{F(d_{i_{1}}, d_{i_{2}})}{f(1, np, np)-f(2, np, np)}-p$. The $ij$-entry of $\bar{A}_{1}$ is\\
\begin{displaymath}
\bar{A}_{1}(i, j)=
\begin{cases}
\bar{F}(d_{i_{1}}, d_{i_{2}}),&\quad \mbox{with probability $p$,}\\
-p,&\quad \mbox{with probability $1-p$}.
\end{cases}
\end{displaymath}
Then,
$$\begin{array}{lll}
\mathbf{E}M_{k}(n^{-1/2}\bar{A}_{1})
&=&\mathbf{E}\frac{1}{n}tr(n^{-1/2}\bar{A}_{1})^{k}\\
&=&n^{-1-k/2}\sum\limits_{i_{1}=1}^{n}\cdots\sum\limits_{i_{k}=1}^{n}\mathbf{E}\{\bar{A}_{1}(i_{1}, i_{2})\cdot \bar{A}_{1}(i_{2}, i_{3})\cdots \bar{A}_{1}(i_{k}, i_{1})\}\\
&=&n^{-1-k/2}\sum\limits_{i_{1}=1}^{n}\cdots\sum\limits_{i_{k}=1}^{n}\mathbf{E}\{a_{i_{1}i_{2}}a_{i_{2}i_{3}}\cdots a_{i_{k}i_{1}}\frac{\bar{A}_{1}(i_{1}, i_{2})}{a_{i_{1}i_{2}}}\cdot
\frac{\bar{A}_{1}(i_{2}, i_{3})}{a_{i_{2}i_{3}}}\cdots \frac{\bar{A}_{1}(i_{k}, i_{1})}{a_{i_{k}i_{1}}}\}.
\end{array}$$

Here, $i.i.d.$ random variables $a_{i_{l}i_{l+1}}=-p$ with probability $1-p$ and $a_{i_{l}i_{l+1}}=1-p$ with probability $p$. Denote by ${WA}_{k}$ the set of $k$-step closed walks on $\{1, 2, \ldots, n\}$. Then
\begin{displaymath}
\begin{split}
\mathbf{E}M_{k}(n^{-1/2}\bar{A}_{1})=&n^{-1-k/2}\sum\limits_{w\in {WA}_{k}}\mathbf{E}\{a_{i_{1}i_{2}}a_{i_{2}i_{3}}\cdots a_{i_{k}i_{1}}\frac{\bar{A}_{1}(i_{1}, i_{2})}{a_{i_{1}i_{2}}}\cdot
\frac{\bar{A}_{1}(i_{2}, i_{3})}{a_{i_{2}i_{3}}}\cdots \frac{\bar{A}_{1}(i_{k}, i_{1})}{a_{i_{k}i_{1}}}\}\\
&=n^{-1-k/2}\sum\limits_{w\in {WA}_{k}}\sum\{a_{i_{1}i_{2}}a_{i_{2}i_{3}}\cdots a_{i_{k}i_{1}}\frac{\bar{A}_{1}(i_{1}, i_{2})}{a_{i_{1}i_{2}}}\cdot
\frac{\bar{A}_{1}(i_{2}, i_{3})}{a_{i_{2}i_{3}}}\cdots \frac{\bar{A}_{1}(i_{k}, i_{1})}{a_{i_{k}i_{1}}}\}\cdot\\
&\mathbf{P}(d_{i_{1}}, d_{i_{2}}, \cdots d_{i_{k}}, a_{i_{1}i_{2}}, a_{i_{2}i_{3}}, \cdots a_{i_{k}i_{1}})\\
&=n^{-1-k/2}\sum\limits_{w\in {WA}_{k}}\sum\{a_{i_{1}i_{2}}a_{i_{2}i_{3}}\cdots a_{i_{k}i_{1}}\frac{\bar{A}_{1}(i_{1}, i_{2})}{a_{i_{1}i_{2}}}\cdot
\frac{\bar{A}_{1}(i_{2}, i_{3})}{a_{i_{2}i_{3}}}\cdots \frac{\bar{A}_{1}(i_{k}, i_{1})}{a_{i_{k}i_{1}}}\}\\
&\cdot\mathbf{P}(d_{i_{1}}, d_{i_{2}}, \cdots d_{i_{k}}|a_{i_{1}i_{2}}, a_{i_{2}i_{3}}, \cdots a_{i_{k}i_{1}})\cdot\mathbf{P}(a_{i_{1}i_{2}}, a_{i_{2}i_{3}}, \cdots a_{i_{k}i_{1}}).
\end{split}
\end{displaymath}
\begin{lemma} (see \cite{boll}) \label{lem2.2} Let $\varepsilon>0$ be fixed, $\varepsilon n^{-3/2}\leq p\leq 1-\varepsilon n^{-3/2}$. Let $q=q(n)$ be a natural number and set\\
\begin{displaymath}
\mu_{q}=nB(q;n-1, p)\quad and\quad\nu_{q}=n\{1-B(q+1;n-1, p)\},
\end{displaymath}
where\\
\begin{displaymath}
B(l;m, p)=\sum\limits_{j\geq l}b(j;m, p)
\end{displaymath}
in which $b(j;m, p)=\binom{m}{j}p^{j}(1-p)^{m-j}$ is subject to the binomial distribution. For a random graph $G\in\mathcal{G}_{n, p}$, denote by $Y_{q}(G)$ the number of vertices of degrees at least $q$ and $Z_{q}(G)$ the number of vertices of degrees at most $q$. Then\\
\begin{displaymath}
(i) \ \  if\,\,\,\mu_{q}\rightarrow 0,\,P(Y_{q}=0)\rightarrow 0;\,
(ii) \ \ if\,\,\,\nu_{q}\rightarrow 0,\,P(Z_{q}=0)\rightarrow 0.
\end{displaymath}
\end{lemma}

Remember that $p\in (0,1)$ is a constant. It is not difficult to check that the minimum and maximum degrees $\delta$ and $\Delta$
of a random graph $G_p$ on $n$ vertices satisfy that
\begin{equation}\label{degree}
np-n^{\frac{3}{4}}<\delta(G_p)\leq\Delta(G_p)<np+n^{\frac{3}{4}},\mbox{ a.s.}
\end{equation}
($(i)$ and $(ii)$ hold by Chernoff's Inequality.) So, we just need to deal with these graphs, in which all vertex degree fall in the interval $(np-n^{\frac{3}{4}},np+n^{\frac{3}{4}})$.
In this case, for all values of $\bar{A}_{1}(i, j)$ in $\bar{A}_{1}$,  we have $\frac{\bar{A}_{1}(i, j)}{a_{i, j}}=1+o(1).$ Then, we can get that

\begin{align}
&n^{-1-k/2}\sum\limits_{w\in {WA}_{k}}\sum_{np-n^{\frac{3}{4}}\leq d_{i_{1}}, d_{i_{2}}, \cdots, d_{i_{k}}\leq np+n^{\frac{3}{4}}}\sum_{a_{i_{1}i_{2}},\cdots,a_{i_{k}i_{1}}}\nonumber\\
&a_{i_{1}i_{2}}a_{i_{2}i_{3}}\cdots a_{i_{k}i_{1}}\frac{\bar{A}_{1}(i_{1}, i_{2})}{a_{i_{1}i_{2}}}\cdot
\frac{\bar{A}_{1}(i_{2}, i_{3})}{a_{i_{2}i_{3}}}\cdots \frac{\bar{A}_{1}(i_{k}, i_{1})}{a_{i_{k}i_{1}}}\cdot\nonumber\\
&\mathbf{P}(d_{i_{1}}, d_{i_{2}}, \cdots d_{i_{k}}|a_{i_{1}i_{2}}, a_{i_{2}i_{3}}, \cdots a_{i_{k}i_{1}})\cdot\mathbf{P}(a_{i_{1}i_{2}}, a_{i_{2}i_{3}}, \cdots a_{i_{k}i_{1}})\nonumber\\
&=n^{-1-k/2}\sum\limits_{w\in {WA}_{k}}\sum_{a_{i_{1}i_{2}},\cdots,a_{i_{k}i_{1}}}a_{i_{1}i_{2}}a_{i_{2}i_{3}}\cdots a_{i_{k}i_{1}}\cdot\nonumber\\
&\mathbf{P}(a_{i_{1}i_{2}}, a_{i_{2}i_{3}}, \cdots a_{i_{k}i_{1}})\sum_{np-n^{\frac{3}{4}}\leq d_{i_{1}}, d_{i_{2}}, \cdots, d_{i_{|V_{w}|}}\leq np+n^{\frac{3}{4}}}(1+o(1))\mathbf{P}(d_{i_{1}}, d_{i_{2}}, \cdots d_{i_{k}}|a_{i_{1}i_{2}}, a_{i_{2}i_{3}}, \cdots a_{i_{k}i_{1}})\nonumber\\
&=n^{-1-k/2}\sum\limits_{w\in {WA}_{k}}\sum_{a_{i_{1}i_{2}},\cdots,a_{i_{k}i_{1}}}(1+o(1))a_{i_{1}i_{2}}a_{i_{2}i_{3}}\cdots a_{i_{k}i_{1}}\mathbf{P}(a_{i_{1}i_{2}}, a_{i_{2}i_{3}}, \cdots a_{i_{k}i_{1}})\nonumber\\
&=n^{-1-k/2}\sum\limits_{w\in {WA}_{k}}(1+o(1))\mathbf{E}(a_{i_{1}i_{2}}a_{i_{2}i_{3}}\cdots a_{i_{k}i_{1}})\nonumber\\
&=n^{-1-k/2}\sum\limits_{w\in {WA}_{k}}\mathbf{E}(a_{i_{1}i_{2}}a_{i_{2}i_{3}}\cdots a_{i_{k}i_{1}})+o(1)\\
&=\mathbf{E}M_{k}(n^{-1/2}\bar{A_{p}})+o(1)\nonumber.
\end{align}

Since $\mathbf{Var}M_{k}(n^{-1/2}\bar{A}_{1})=\mathbf{E}M_{k}^{2}(n^{-1/2}\bar{A}_{1})-\{\mathbf{E}M_{k}(n^{-1/2}\bar{A}_{1})\}^{2}$, repeating the process above, similarly we can get\\
\begin{displaymath}
\mathbf{Var}(M_{k}(n^{-1/2}\bar{A}_{1}))=\mathbf{Var}(M_{k}(n^{-1/2}\bar{A}_{p}))+o(1)\quad a.s.
\end{displaymath}
It was proved in \cite{bai, wig1} that\\
\begin{displaymath}
\mathbf{Var}(M_{k}(n^{-1/2}\bar{A}_{p}))=O(n^{-2})\quad for \ all \ k>0.
\end{displaymath}
In other words, the fluctuations of $M_{k}(n^{-1/2}\bar{A_{1}})$ and $M_{k}(n^{-1/2}\bar{A}_{p})$ vanish when $n\rightarrow\infty$. Therefore,  $\mathbf{E}M_{k}(n^{-1/2}\bar{A_{1}})$ and $\mathbf{E}M_{k}(n^{-1/2}\bar{A}_{p})$ can be approximately regarded as $M_{k}(n^{-1/2}\bar{A_{1}})$ and $M_{k}(n^{-1/2}\bar{A}_{p})$. That is, $\lim\limits_{n\rightarrow\infty}\mathbf{E}M_{k}(n^{-1/2}\bar{A_{1}})=\int x^{k}\,d\Phi(x)$\quad a.s. for any $k$. Using Lemma \ref{lem1.5}, the limiting distribution of eigenvalues of $n^{-1/2}\bar{A_{1}}$ is the same as $n^{-1/2}\bar{A}_{p}$. According to Theorem \ref{thm1}, \\
\begin{displaymath}
\lim\limits_{n\rightarrow\infty}\Phi_{n^{-1/2}\bar{A}_{1}(x)}=\Phi(x)\quad a.s.
\end{displaymath}
As a consequence, \\
\begin{align}
\mathscr{E}(\bar{A_{1}})&= n^{\frac{3}{2}}\int |x|\,d\Phi_{n^{-1/2}\bar{A}_{1}(x)}\nonumber\\
&= n^{\frac{3}{2}}(\int |x|\,d\Phi(x)+o(1))\nonumber\\
&=(\frac{8}{3\pi}\sqrt{p(1-p)}+o(1))\cdot n^{\frac{3}{2}}\nonumber
\end{align}

\begin{remark}
The proof of equation (2) is from \cite{bai}. We write\\
\begin{align}
&n^{-1-k/2}\sum\limits_{w\in {WA}_{k}}(1+o(1))\mathbf{E}(a_{i_{1}i_{2}}a_{i_{2}i_{3}}\cdots a_{i_{k}i_{1}})=n^{-1-k/2}\sum\limits_{w\in {WA}_{k}}\mathbf{E}(a_{i_{1}i_{2}}a_{i_{2}i_{3}}\cdots a_{i_{k}i_{1}})\nonumber\\
&+o(n^{-1-k/2}\sum\limits_{w\in {WA}_{k}}|\mathbf{E}(a_{i_{1}i_{2}}a_{i_{2}i_{3}}\cdots a_{i_{k}i_{1}})|).\nonumber
\end{align}
In the proof of Theorem 2.1 of \cite{bai}, Bai studied the sum $n^{-1-k/2}\sum\limits_{w\in {WA}_{k}}\mathbf{E}(a_{i_{1}i_{2}}a_{i_{2}i_{3}}\cdots a_{i_{k}i_{1}})$ and showed that: if $k=2s+1$, the number of non-vanishing terms is no more than $n^{s+1}$; if $k=2s$, the number of negative terms is no more than $n^s$. According to Lemma \ref{lem2.1}, the sum in the brackets is bounded as $n\rightarrow\infty$, and thus equation (2) is obtained. The variance case can be verified similarly.
\end{remark}
\begin{remark}
Although $|x|$ does not have a compact support set, it is right that
$$\int |x|\,d\Phi_{n^{-1/2}\bar{A}_{1}(x)}\rightarrow\int |x|\,d\Phi(x).$$ One can check it by dividing $\mathbf{R^{1}}$ into two parts: an interval $\emph{I}$ with $[-2, 2]\subseteq\emph{I}$ and the remaining unbounded segments.
\end{remark}
\begin{lemma} \label{fan}
(Ky Fan \cite{fan}) Let $X, Y$ and $Z$ be real symmetric matrices of order $n$ and $X+Y=Z$. Then\\
\begin{displaymath}
\mathscr{E}(X)+\mathscr{E}(Y)\geq\mathscr{E}(Z).
\end{displaymath}
\end{lemma}
According to Lemma \ref{fan}, $$ \mathscr{E}(\bar{A_1})-\mathscr{E}(-J+I)\leq \mathscr{E}(\frac{A_{1}}{f(1, np, np)-f(2, np, np)})\leq \mathscr{E}(\bar{A_1})+\mathscr{E}(J-I) . $$

We know  that $\mathscr{E}(J-I)=2(n-1)$. Consequently, $$\mathscr{E}(A_{1})=|f(1, np, np)-f(2, np, np)|(\frac{8}{3\pi}\sqrt{p(1-p)}+o(1))\cdot n^{\frac{3}{2}}.$$
Now, we finished the estimation of the energy of  $A_{1}$. In next subsection, we focus on the other matrix $A_2$ to get our final result.

\subsection{Estimation of $\mathscr{E}(W_f(G_p))$}

We write $A_{2}=A_{2}^{\prime}+A_{2}^{\prime\prime}$, where the $ij$-entries $A_{2}^{\prime}(i,j)=f(2, np, np)$, $A_{2}^{\prime\prime}(i, j)=f(2,d_i,d_j)-f(2, np, np)$ with $i\neq j$, and $A_{2}^{\prime}(i,j)=0, A_{2}^{\prime\prime}(i, j)=0$ with $i=j$. We know that $$\mathscr{E}(A_{2}^{\prime})=2(n-1)|f(2, np, np)|.$$
Denote the eigenvalues of $A_{2}^{\prime\prime}$ by $\theta_{1}, \ldots,\theta_{n}$ . It is easy to see that $$\theta_{1}^{2}+\cdots+\theta_{n}^{2}=trA_{2}^{\prime\prime 2}=\sum\limits_{i,j}A_{2}^{\prime\prime 2}(i, j)\leq n(n-1)\cdot\max_{i\neq j}(f(2,d_i,d_j)-f(2,np,np))^2.$$
Recall that formula (\ref{degree}) holds for almost all graphs, and then $ \max_{i\neq j}|f(2,d_i,d_j)-f(2,np,np)|=o(1)|f(2,np,np)|.$
By Cauchy-Schwarz Inequality,
\begin{displaymath}
\mathscr{E}(A_{2}^{\prime\prime})=|\theta_{1}|+\cdots +|\theta_{n}|\leq\sqrt{n}\cdot\sqrt{n(n-1)}\cdot o(1)|f(2, np, np)|=o(n^{\frac{3}{2}})|f(2, np, np)|.
\end{displaymath}
Using Lemma \ref{fan}, $\mathscr{E}(A_{2})=o(n^{\frac{3}{2}})|f(2, np, np)|$  for almost all graphs.

Combining $\mathscr{E}(A_{1})$ and $\mathscr{E}(A_{2})$, and applying Lemma \ref{fan} again, the proof of Theorem \ref{main} is thus finished.

Moreover, the special case $f(1, np, np)=f(2, np, np)$ can be solved by the same method above, and then there is no need to partition $W_f(G_p)$ into two matrices, and instead, we can get the estimation $\mathscr{E}(W_f(G_p))$ directly.

\section[Applications for matrices with distance-based and degree-distance-based weights of chemical use]{Applications for matrices with distance-based and degree-distance based weights of chemical use}

Topological indices in chemistry are used to represent structural properties of molecular graphs. Each index maps a molecular graph into a single number, obtained by summing up the weights of all pairs of vertices in a molecular graph. If we use a matrix to represent the structure of a molecular graph with weights separately on its pairs of vertices,
it will completely keep the structural information of the graph, i.e., a matrix keeps much more structural information than an index. So, further study on the algebraic properties of
these structural matrices should be made in the future.

Now we can get the asymptotic values of energies for various kinds of matrices with special degree-distance-based weights. Here we only deal with matrices with distance-based weights and matrices with general degree-distance-based weights, since many examples of matrices with only degree-based weights were given in \cite{lls}. We also recommend \cite{furt} for more such weight functions from indices.

The first that should be mentioned is the distance matrix of a graph $G$, coming from the famous Wiener index. The distance matrix $dist(G)$ is the matrix of $G$ with the $ij$-entry weighted by the distance between vertices $i$ and $j$, that is, $f(D(i, j), d_{i}, d_{j})=D(i, j)$. Applying Theorem \ref{main}, we can get \\
\begin{corollary}
\begin{displaymath}
\mathscr{E}(dist(G))=(\frac{8}{3\pi}\sqrt{p(1-p)}+o(1))n^{3/2}\quad a.s.
\end{displaymath}
\end{corollary}
This can be verified by \cite{dll2}.\\
\\
The Harary matrix \cite{plav, ivanc} $RD(G)$ is defined as\\
\begin{displaymath}
RD(i, j)=
\begin{cases}
1/D(i, j),&\quad \mbox{if $i\neq j$}\\
0,&\quad \mbox{if $i=j$}
\end{cases}
\end{displaymath}
and the energy of $RD(G)$ is called Harary energy, denoted by $HE(G)$.
From Theorem \ref{main}, we can get\\
\begin{corollary}
\begin{displaymath}
HE(G)=(\frac{4}{3\pi}\sqrt{p(1-p)}+o(1))n^{3/2}\quad a.s.
\end{displaymath}
\end{corollary}

The hyper-Wiener index, introduced in 1993 by Randi$\acute{c}$ in [18],  is defined as $\frac{1}{2}\sum\limits_{u, v\in V(G)}(D(u, v)+D^{2}(u, v))$. Similarly, we can define the hyper-Wiener matrix $HW(G)$ with $ij$-entry as follows:\\
\begin{displaymath}
HW(G)(i, j)=\frac{1}{2}(D(i, j)+D^{2}(i, j)).
\end{displaymath}
From Theorem \ref{main}, we can get
\begin{corollary}
\begin{displaymath}
\mathscr{E}(HW(G))=(\frac{16}{3\pi}\sqrt{p(1-p)}+o(1))n^{3/2}\quad a.s.
\end{displaymath}
\end{corollary}
Denote the diameter of $G$ by $D(G)$. The reciprocal complementary Wiener index, introduced in 2000 by Ivanciuc et al. in \cite{iib}, is defined as $\sum\limits_{u, v\in V(G)}\frac{1}{D(G)+1-D(u, v)}$. We then can define the reciprocal complementary Wiener matrix $RCW(G)$ with $ij$-entry as follows:\\
\begin{displaymath}
RCW(G)(i, j)=
\begin{cases}
\frac{1}{D(G)+1-D(i, j)}, &\quad \mbox{if $i\neq j$}\\
0,&\quad \mbox{if $i=j$}.
\end{cases}
\end{displaymath}
From Theorem \ref{main}, we can get
\begin{corollary}
\begin{displaymath}
\mathscr{E}(RCW(G))=(\frac{4}{3\pi}\sqrt{p(1-p)}+o(1))n^{3/2}\quad a.s.
\end{displaymath}
\end{corollary}
The reverse Wiener matrix \cite{van} $RW(G)$ is defined as\\
\begin{displaymath}
RW(i, j)=
\begin{cases}
D(G)-D(i, j),&\quad \mbox{if $i\neq j$}\\
0,&\quad \mbox{if $i=j$}.
\end{cases}
\end{displaymath}
From Theorem \ref{main}, we can get
\begin{corollary}
\begin{displaymath}
\mathscr{E}(RW(G))=(\frac{8}{3\pi}\sqrt{p(1-p)}+o(1))n^{3/2}\quad a.s.
\end{displaymath}
\end{corollary}

As one can see, all the above matrices are defined purely distance-based. Next we show some degree-distance-based ones. The first one is from the  topological index in \cite{dob} defined as $DD=\sum\limits_{u\neq v}(d(u)+d(v))D(u, v)$. We then define the degree-distance-based matrix $DD(G)$ with $ij$-entry as  \\
\begin{displaymath}
DD(G)(i, j)=(d_{i}+d_{j})D(i, j).
\end{displaymath}
From Theorem \ref{main}, we can get
\begin{corollary}
\begin{displaymath}
\mathscr{E}(DD(G))=(\frac{16}{3\pi}p\sqrt{p(1-p)}+o(1))n^{5/2}\quad a.s.
\end{displaymath}
\end{corollary}
The Gutman index \cite{gut2} is defined as $Gut=\sum\limits_{u\neq v}d(u)d(v)D(u, v)$. We then define the Gutman matrix $Gut(G)$ with $ij$-entry as\\
\begin{displaymath}
Gut(G)(i, j)=d_{i}d_{j}D(i, j).
\end{displaymath}
From Theorem \ref{main}, we can get
\begin{corollary}
\begin{displaymath}
\mathscr{E}(Gut(G))=(\frac{8}{3\pi}p^{2}\sqrt{p(1-p)}+o(1))n^{7/2}\quad a.s.
\end{displaymath}
\end{corollary}

The additively weighted Harary index and multiplicatively weighted Harary index were introduced in \cite{ali, hua}, which are defined respectively as follows:\\
\begin{displaymath}
H_{A}=\sum\limits_{u\neq v}\frac{d(u)+d(v)}{D(u, v)},
\end{displaymath}
\begin{displaymath}
H_{M}=\sum\limits_{u\neq v}\frac{d(u)d(v)}{D(u, v)}.
\end{displaymath}
We then define the additively weighted Harary matrix and multiplicatively weighted Harary matrix with $ij$-entry, respectively, as follows:\\
\begin{displaymath}
H_{A}(G)(i, j)=
\begin{cases}
\frac{d_{i}+d_{j}}{D(i, j)},&\quad \mbox{if $i\neq j$}\\
0,&\quad \mbox{if $i=j$}.
\end{cases}
\end{displaymath}
\begin{displaymath}
H_{M}(G)(i, j)=
\begin{cases}
\frac{d_{i}d_{j}}{D(i, j)},&\quad \mbox{if $i\neq j$}\\
0,&\quad \mbox{if $i=j$}.
\end{cases}
\end{displaymath}
From Theorem \ref{main}, we can get
\begin{corollary}
\begin{displaymath}
\mathscr{E}(H_{A}(G))=(\frac{8}{3\pi}p\sqrt{p(1-p)}+o(1))n^{5/2}\quad a.s.
\end{displaymath}
\begin{displaymath}
\mathscr{E}(H_{M}(G))=(\frac{4}{3\pi}p^{2}\sqrt{p(1-p)}+o(1))n^{7/2}\quad a.s.
\end{displaymath}
\end{corollary}

There are variety of topological indices of chemical use, which are
degree-based only, distance-based only, and mixed with both degree and distance based. From them we can get the corresponding weighted matrices. Then from our result Theorem \ref{main} we can get the asymptotic values of energies of all these corresponding random weighted graphs.

\end{document}